
\documentstyle{amsppt}
\magnification=1200
\def\Q{\Bbb Q}
\def\Z{\Bbb Z}
\def\O{\Cal O}
\def\K{\Bbb K}
\def\N{\bold N}
\def\C{\Bbb C}
\def\R{\Bbb R}
\def\F{\Bbb F}
\def\A{A}
\def\W{W}
\def\nth{^{th}}
\def\T{T}
\def\I{I}
\def\GL{GL}
\def\SL{SL}
\def\pigam{\left(2\pi\right)^{1-3w}
   \zeta\left(3w-\frac{4}{3}\right)
   \zeta\left(3w-\frac{2}{3}\right) \Gamma \left(
   \frac{3w}{2} - \frac{2}{3} \right) \Gamma \left( \frac{3w}{2} -
   \frac{1}{3} \right) }

\def\p{\rho}
\def\tp{\tilde{\rho}}
\topmatter
\title Automorphic Forms \\
and Cubic Twists of Elliptic Curves
\endtitle
\rightheadtext{AUTOMORPHIC FORMS AND CUBIC TWISTS}
\author
Daniel Lieman
\endauthor
\address
Columbia University Math Department, New York, NY  10027
\endaddress
\email
lieman\@math.columbia.edu
\endemail
\thanks
The author wishes to thank Jeffrey Hoffstein and Daniel Bump
for much patient advice
and consultation, as well as Dorian Goldfeld,
Morris Newman, Michael Rosen, and Don
Zagier for helpful comments and referrals to relevant literature.  The author
would also like to thank M.S.R.I., Berkeley, for its gracious hospitality
during the completion of some of the
research surveyed in this paper; during that time, he was partially
supported by NSF grant DMS-9022140.  The author was partially supported by
NSF grant DMS-9400882 during final preparation of this article.
This paper is an exposition of a lecture presented to the New York
Number Theory Seminar in the fall of 1993.  The author would like to
thank D. Chudnovsky for the invitation to speak.
\endthanks
\endtopmatter
\document
\head
0. Introduction
\endhead

One of the most classical problems in number theory is that of determining
whether a given rational integer is the sum of two cubes of rational numbers.
This ``sum of two cubes'' problem has been attacked from a variety of
both classical and modern viewpoints; it would be nearly impossible to
list here all of the various approaches taken and results obtained.    An
extensive compilation of the older history of the problem is given
in Dickson~\cite{D}.

This paper is an exposition of one of the connections between the curve
$$E_D: x^3 + y^3 = D \tag{0.1}$$
and metaplectic forms.  The approach we shall consider has the potential of
providing collective
information ``on average'' about the curves $E_D$.  It will not
provide any information about $E_D$ for a particular $D$; for these types of
results, the best modern information seems to be from Elkies, and from
Rodriguez Villegas and Zagier.  Here, we will describe the ingredients
necessary
to obtain two theorems concerning the curves $E_D$.

\proclaim{Theorem 1} There are infinitely many cube--free $D$ such that
$E_D$ has no rational solutions.
\endproclaim

\proclaim{Theorem 2} Fix a prime $p \ne 3$, and any congruence class
$c$ modulo $p$.  Then there are infinitely many cube--free $D$ congruent
to $c$ modulo $p$ such that $E_D$ has no rational solutions.
\endproclaim

\demo{Remark}
Theorem 1 is in fact weaker than a classical result of Sylvester, while
Theorem 2 is only slightly stronger than Sylvester's results combined
with the Dirichlet theorem on primes in arithmetic progressions.  Our
main results here are Theorem 3, relating the L--series to $E_D$ to
a certain metaplectic form, and the machinery of Theorem 4 and the
subsequent discussion, which allow
one to obtain analytic information about these L--series.
\enddemo

There are many connections between $E_D$ and automorphic forms.  The curves
$E_D$ are elliptic curves, and there is the famous conjecture that they are
related to a form of weight~2.  In addition,
work of Waldspurger and the Shimura
correspondence give the existence
of a form $f$ of weight~$3/2$ which is related to the values of the L--series
of $E_D$.  Indeed, Nekovar~\cite{N} has explicitly identified the
corresponding form.  (It it perhaps interesting to note that this form is
a product of Dedekind $\eta$~functions, and that the formula for the
value of the L--series of $E_D$ at the center of the critical strip
given by Rodriguez Villegas and Zagier also depends on $\eta$ functions.
It is not clear whether this is more than coincidence.)

The connection between the L--series of $E_D$ and metaplectic forms
is simultaneously quite explicit and not well--understood.  That is, we
will see that the L--series arise as the Whittaker--Fourier coefficients
of a certain metaplectic Eisenstein series on the cubic cover of $GL(3)$,
and perhaps we will even make sense of $why$ one might expect this to happen;
we will not, however, be able to satisfactorily explain this occurrence.
There is no known analog of Waldspurger's result in this case.

It is clear that this is not an isolated occurrence of an L--series in
the coefficients of a metaplectic form.
Recent works of
Bump and Hoffstein~\cite{BH},
Bump, Friedberg and Hoffstein~\cite{BHF},~\cite{BHF2}
and Goldfeld, Hoffstein, and Patterson~\cite{GHP} have discovered
metaplectic Eisenstein series with Whittaker--Fourier coefficients
which are (essentially) the Hecke L--series of the cubic residue symbol,
the L--series of quadratic twists of elliptic curves with complex
multiplication, and the L--series of quadratic twists of cuspidal
newforms for the group $\Gamma_0(M)$.    For all the cases mentioned above,
one is able to apply various analytic methods to the appropriate metaplectic
form and obtain interesting
number--theoretic results.

This paper will survey the construction of the form with coefficients
which are the L--series of cubic twists of elliptic curves with complex
multiplication, the extraction of information about the curves from this
form, and how much this process may be generalized using a similar
construction to obtain L--series of higher order twists, or Hecke L--series
of higher order residue symbols. It is a summary of the results
of \cite{L},\cite{L2}, and \cite{BL}, in a manner less terse than
that usually found in research articles.

We should mention that we often ignore the question of attribution,
especially for fundamental work, in what follows.  We take for granted
several highly insightful observations and constructions which are
crucial to the general theory of metaplectic forms and the success
in carrying out investigations within this framework.  For a
more detailed and far--reaching survey of the theory of metaplectic
forms, consult the beautiful and quite readable
article of Hoffstein~\cite{H}.

\head
1. Cubic twists of elliptic curves
\endhead

We are interested in the family of curves
$$x^3 + y^3 = D. \tag{1.1}$$
As mentioned above, these are elliptic curves (with Weierstrass form
$y^2 = x^3 - 432D^2$), and the rich machinery of elliptic curves may
be brought to bear.  To each curve there is associated an L--series
$L(E_D,s)$, and the following are well-known.

\proclaim{Theorem (Mordell--Weil)} The set of rational solutions to
$E_D,$ together with a ``point at infinity,'' form a finitely generated
abelian group $E_D(\Q)$, under a certain geometric group law.
\endproclaim

\proclaim{Theorem (Coates--Wiles)}
If $L(E_D,1) \ne 0,$ then $E_D$ has only finitely
many rational points.
\endproclaim
\proclaim{Proposition} For $|D| \ge 3$, the torsion subgroup
$E_D(\Q)_{\text{tors}} \subseteq E_D(\Q)$ is trivial.
\endproclaim

The Mordell--Weil Theorem asserts that
$$E_D(\Q) \cong \Z^{r_D} \oplus E_D(\Q)_{\text{tors}}.$$
Combined with the proposition, this shows that for all but a few $D$
(in particular, $D=0, \pm 1, \pm 2$), one has
$$E_D(\Q) \cong \Z^{r_D}. \tag{1.2}$$
Now applying the Theorem of Coates--Wiles to (1.2), we obtain the
statement
$$L(E_D,1) \ne 0 \,\text{ implies }\, r_D = 0. \tag{1.3}$$
It is the behavior of $r_D$ that we want to study, and we shall do so using
the relationship (1.3).

We begin by setting some notation.  Let $\K$ denote, for the rest of the
paper, $\Q(\sqrt{-3})$, and let $\O$ denote the ring of integers in $\K$.
If $a,b \in \O, b \equiv 1\,(3)$, we write $\left(\dfrac{a}{b}\right)_3$
for the cubic residue symbol of $a$ mod $b$.  The following
proposition follows from the
exposition in Ireland and Rosen~\cite{IR}, together with the Weierstrass
form of $E_D$, mentioned above.
\proclaim{Proposition} The L--series of the elliptic curve $E_D$ is given
by $$L(E_D,s) = \sum \Sb a \in \O \\ a \equiv 1 \, (3) \endSb
\left(\frac{D}{a}\right)_3 \frac{|a|}{a}\,\N(a)^{-s + \frac{1}{2}}. \tag{1.4}$$
\endproclaim

\demo{Remark}
Each of the curves $E_D$ is a cubic twist of the curve $E_1$, and so the
proposition says that the L--series is a natural object, that is, that the
L--series of ``$E_1$ twisted by $D$'' is just the ``L--series of $E_1$''
twisted by the cubic residue symbol $\left(\frac{D}{\cdot}\right)_3$.
On the other hand, one can also interpret the Proposition as saying that
the L--series of $E_D$ is just the Hecke L--series of the cubic
residue symbol $\left(\frac{D}{\cdot}\right)_3$ twisted by the
$grossencharacter$
$\frac{|\cdot|}{\cdot}$.  This latter viewpoint is the basis of our
arguments below.
\enddemo

\demo{Remark}
The family of curves we are studying is a particularly interesting one,
not merely for its long history, but also because of recent results.  Zagier
and Kramarz~\cite{ZK} have computed the analytic rank of $E_D$ for $D$
cube--free
and less than $70,000$ and have found that the distribution of curves with
analytic rank $\ge 2$ is quite frequent ($\thicksim 25\%$).  More impressively,
they found that, writing $X^*$ for the number of cube--free integers less
than $X$, the average
$$\frac{1}{X^*} \cdot \#\{D < X \, | \, D \text{ cube--free, with analytic rank
of
$E_D \ge 2$}\}$$ seemed to
be independent of $X$ ($X$ ranging up to $70,000$).  This rather uniform
distribution of curves of high analytic rank is somewhat surprising, and it
is an interesting question to see to what extent, if at all,
this is reflected in the average
values of the L--series of these curves.
\enddemo

\demo{Remark}
The average value of the L--series of $E_D$ as $D$ varies is the
subject of two different conjectures.
Zagier~\cite{Z-K}
has pointed out that the series  $$L_{av}(s) = \sum_n \frac{b_n}{n^s}$$  formed
by setting $b_n$ to be the average of the $n\nth$ coefficient
of the L--series of $E_D$ as $D$ varies
has an analytic continuation with a finite
value at s=1, and that this makes it plausible that
the numbers  $L(E_D,1)$  as $D$  varies have a well--defined average value.
Goldfeld and Viola~\cite{GV}, on the other hand, have given (based on
heuristic arguments) a very
general conjecture, the specialization of which to this problem is
$$\sum \Sb p < X \\ p \text{ prime } \endSb
L(E_p,1) \,L(E_{p^2},1) \thicksim c \cdot X  \text{ as } X \to \infty.$$
(Goldfeld and Viola also explicitly calculate the expected value of
the constant $c$).  We will discuss average values further in section~5.
\enddemo

\head
2. Overview of the approach
\endhead

We begin by recalling how analytic information
about the behavior of a Dirichlet series can yield information about
the behavior of the coefficients.  By a Dirichlet series, we mean
a function of the form
$$f(s) = \sum_{n=1}^\infty \frac{a_n}{n^s}$$
satisfying
$$\align
(1) \qquad&
\text{$f(s)$ converges in a half-plane (i.e. for $re(s)$ sufficiently large)}
\\
(2) \qquad&
\text{$f(s)$ has a meromorphic continuation to the entire
complex plane.}
\endalign$$
The prototype to think of, of course, is the Riemann $\zeta-$function,
$$\zeta(s) = \sum_{n=1}^\infty \frac{1}{n^s} = \prod_p (1-p^{-s})^{-1}$$
which satisfies the functional equation
$$\zeta^*(s) = \Gamma\left(\frac{s}{2}\right) \pi^{-\frac{s}{2}} \zeta(s) =
\zeta^*(1-s).$$

There are two ways a Dirichlet series, such as the $\zeta-$function, can
yield information about its coefficients.

\proclaim{The easy method (Obvious)} Suppose $f(s)$ is a Dirichlet series, as
above,
and that $f(s)$ has a pole at $s=s_0$.  Then infinitely many of the
coefficients $a_n$ are non-zero.
\endproclaim

\proclaim{The hard method (Ikehara's Tauberian Theorem)} Suppose $f(s)$ is a
Dirichlet series which converges for $re(s) > 1$, and has a
pole of order $1$ at $s=1$ with residue~$\alpha$.  Further suppose
that $f(s)$ is analytic on the line $re(s)~=~1$ except for the pole
at $s=1$, and that
each $a_n \ge 0$.  Then
$$\sum_{n<X} a_n \thicksim \alpha X \qquad
\text{as $X \to \infty$}.$$
\endproclaim

\demo{Remark}
Tauberian theorems are of course much more general than
that given above.  We use the simplest specialization to illustrate the general
shape of such a theorem, and because the $\zeta-$function satisfies the
hypothesis, allowing the reader to verify the truth of the
method in this one case.
For this paper, we will use only ``the easy method.''   We will discuss
expected applications of Ikehara's Tauberian theorem
to our problem, but without results.
\enddemo

We are now able to outline, in more detail, the contents of this paper.
Our first topic will be the metaplectic cover of $GL(r)$.  We will
define automorphic functions on this group (the metaplectic group),
and see that built into the
automorphy property
is the $n^{th}$ order reciprocity symbol, thus making it plausible that the
predicted L--series should occur in the coefficients of these forms.
We will then construct a particular form on the cubic cover of $GL(3)$,
with the property that its Whittaker--Fourier coefficients are
the L--series of the elliptic curves, in the form given in (1.4).  We will next
use the Rankin--Selberg method to construct a Dirichlet series (more or less)
of the form
$$\sum_{D=1}^\infty \frac{L(E_D,1)}{D^w}$$
and to obtain analytic information about the behavior of this Dirichlet
series in the variable $w$.
Finally, we will discuss generalizations of this work: what is known, and
what is conjectured for L--series of higher order twists, and higher order
residue symbols.

\head
3. Automorphic forms on the metaplectic group
\endhead

There are several ways to describe functions on the metaplectic group.
We take here the most concrete, though perhaps not the most direct; we will
describe a multiplier system, and show that functions transforming with
respect to this multiplier system are in fact functions on the double cover
of $GL(2)$.  We will then discuss how this construction
generalizes to functions on
the $n-$cover of $GL(r)$.  Our construction will
clearly embed the quadratic reciprocity law into the forms we will construct;
it is this key point which explains, in part, why one might expect the
Hecke L--series of quadratic residue symbols (or quadratic
gauss sums, for that matter) to appear in the Fourier
coefficients of a form constructed in this manner.
A basic reference for this section is Hoffstein~\cite{H}.

We begin by defining forms on the double cover of $GL(2)$. The first
step is to construct the Kubota homomorphism.
Let $\lambda = 1+i$, with $i=\sqrt{-1}$, as usual,
and write $\left(\dfrac{a}{b}\right)_2$ for the quadratic residue symbol.
We let $\Gamma(\lambda^3)$ denote the principal congruence subgroup modulo
$\lambda^3$, that is
$$\Gamma(\lambda^3) = \left\{ \pmatrix a & b \\ c & d \endpmatrix \in
SL(2,\Q(i)) \, | \, a \equiv d \equiv 1 \, (\lambda^3), b \equiv c \equiv 0
\, (\lambda^3)\right\}.$$
We define a homomorphism $\kappa: \Gamma(\lambda^3) \to \{\pm1\}$ by
$$\kappa\pmatrix a&b\\c&d\endpmatrix = \cases \left(\dfrac{a}{b}\right)_2
& \quad \text{ if } c \ne 0, \\
\,\,\,1& \quad \text{ if } c = 0. \endcases$$
That this is indeed a homomorphism follows from
quadratic reciprocity.

The standard action of $GL(2,\R)$ on the upper half plane may be generalized
in the following way.  Let $Z~\subset~GL(2,\R)$ and
$K~\subset~GL(2,\R)$ denote the groups of
scalar and orthogonal matrices, respectively; each coset
$GL(2,\R)/Z\cdot K$ has a unique representative of the form $\pmatrix y&x \\
0& 1 \endpmatrix$, where $x,y~\in~\R$, and $y > 0$ (this decomposition
is called the Iwasawa decomposition; each such coset representative is called
the Iwasawa coordinate).  Further, if we let
$GL(2,\R)$ act on the upper half plane in the usual way (i.e.
$\pmatrix a & b \\ c & d \endpmatrix:z~\mapsto~\dfrac{az+b}{cz+d}$), then for
$z = x + iy$, and $\gamma \in GL(2,\R)$, we have $\gamma z = x' + iy'$, where
$\pmatrix y' & x' \\ 0 & 1 \endpmatrix$ is the Iwasawa coordinate of
$\gamma \pmatrix y & x \\ 0 & 1 \endpmatrix$.  That is, the action of
$GL(2,\R)$
on the upper half plane is the same as the action of $GL(2,\R)$ on
$GL(2,\R)/Z\cdot K$ by left multiplication.

This latter action is easy to generalize to larger fields, or to higher
rank groups.  For now, we consider the first generalization.  Let
$Z$ and $K$ now denote, respectively,
the subgroups of scalar and unitary matrices within
$GL(2,\C)$.  We may consider complex--valued functions on $GL(2,\C)/Z \cdot K$
which satisfy
$$f(\gamma \tau) = \kappa(\gamma)f(\tau) \text{ for all }
\gamma \in \Gamma(\lambda^3), \tau \in GL(2,\C)/Z \cdot K. \tag{*}$$

We now explain the correspondence between functions satisfying (*) and
functions on the double cover of $GL(2,\C)$.  We form the double
cover $$\tilde{G}~=~\{(g,\epsilon)\,|\,g\in~GL(2,\C),\epsilon=\pm1\}$$ with
componentwise multiplication
$(g,\epsilon)(g',\epsilon')=(gg',\epsilon\epsilon')$.  We also fix the
subgroups $$\tilde{\Gamma}=\{(g,\kappa(g))\,|\,g\in~\Gamma(\lambda^3)\},
\qquad \tilde{Z}=\{(z,1)\,|\,z~\in~Z\}, %
\qquad \tilde{K}=\{(k,1)\,|\,k~\in~K\}.$$
Consider now a function $\tilde{f}$ on $\tilde{G}/\tilde{Z}\cdot\tilde{K}$
satisfying: (1) $\tilde{f}$ is invariant under the action of $\tilde{\Gamma}$,
and (2) $\tilde{f}((g,\epsilon))~=\epsilon\tilde{f}((g,1)).$  Finally,
define a function $f$ on $GL(2,\C)/Z\cdot K$ by $f(g) = \tilde{f}((g,1)).$
Then we have, for $\gamma~\in~\Gamma(\lambda^3)$,
$$\multline
\qquad \qquad f(\gamma g) = \tilde{f}((\gamma g,1)) = \tilde{f}((\gamma g,
\kappa(\gamma)^2))
= \tilde{f}((\gamma,\kappa(\gamma))(g,\kappa(\gamma))) \\
= \tilde{f}((g,\kappa(\gamma))) = \kappa(\gamma)\tilde{f}((g,1)) =
\kappa(\gamma) f(g).\qquad \qquad \endmultline$$
Thus functions on the double cover of $GL(2,\C)$ which satisfy the conditions
stated above correspond to functions on $GL(2,\C)/Z \cdot K$ satisfying
condition (*).  (Functions on the double cover of $GL(2,\C)$ satisfying the
first condition (regardless of whether they satisfy the automorphy (second)
condition) are called $genuine$.)
$Metaplectic$ $forms$ are functions which satisfy (*) and which also
satisfy a particular differential equation.  This latter condition is
a technical one, and we do not dwell on it here.

We now wish to show how to use a similar construction to find functions
on the $n-$cover of $GL(r,\C)$. It is clear how we would define functions
on the $n-$cover of $GL(2,\C)$; we would use the same construction as
above, except that $\kappa$ would be derived from the $r^{th}$ order
residue symbol.  For higher rank groups, it becomes difficult to
describe the homomorphism $\kappa$ necessary to define forms on covering
group.  Nonetheless, the $n-$cover of $GL(r)$ is defined, and metaplectic
forms (automorphic forms on this covering group) do exist; the
concrete framework we have chosen for our discussion is the wrong setting
to discuss general metaplectic forms (see~\cite{KP} for a discussion
of metaplectic forms on the covers of $GL(r)$).

One case in which $\kappa$ has been explicitly written down, for higher
rank groups, is the case of the cubic cover of $GL(3,\C)$.  Because this
is important for our applications, we review the construction here
(details may be found in~\cite{BH}).  As before, we write $Z$ and $K$ for the
scalar matrices and unitary
matrices of $GL(3,\C)$, respectively; we will wish to define functions
on the space $GL(3,\C)/Z\cdot K$ which satisfy
$f(\gamma\tau)=\kappa(\gamma)f(\tau)$ for $\gamma$ in an appropriate
discrete subgroup $\Gamma,$ as before, and $\kappa:\Gamma\to\mu_3$ an
appropriate homomorphism which takes values in the set of
cube roots of unity.
Such functions, as we have seen, correspond to functions on the cubic
cover of $GL(3,\C)$.  Each coset of $GL(3,\C)/Z\cdot K$ has a unique
Iwasawa representative, as before, of the form
$$\pmatrix y_1y_2 & & \\ & y_1 & \\ & & 1 \endpmatrix
\pmatrix 1 & x_2 & x_3 \\ & 1 & x_1 \\ & & 1 \endpmatrix,$$
where the $y_i$ are positive real numbers, and the $x_i$ are complex
numbers.  We may thus write down our desired functions concretely as
functions of two real and three complex variables; this is precisely
what we will do in the next section.

Turning now to the definition of $\kappa$,
we recall our earlier notation of
$\K$ and $\O$ for $\Q(\sqrt{-3})$ and its ring of integers, respectively.
We write $\Gamma$ for the subgroup of all matrices in $SL(3,\O)$ congruent
to the identity modulo 3.  We need a way of parametrizing matrices
in $\Gamma$; to identify such a system, we first note that that
$GL(3,\C)$ possesses an involution
$$\iota: g \mapsto \pmatrix & & 1 \\ & 1 & \\ 1 & & \endpmatrix \cdot
{}^tg^{-1} \cdot \pmatrix & & 1 \\ & 1 & \\ 1 & & \endpmatrix.$$
The bottom rows of the matrices $g$ and $\iota(g)$ depend only on
the orbit of $g$ in $$GL(3,\C)_\infty~\setminus~GL(3,\C),$$ where
$GL(3,\C)_\infty$ denotes the upper triangular unipotent matrices.
We denote the bottom row of $g$ by the coordinates $A_1$, $B_1$, $C_1$,
and the bottom row of $\iota(g)$ by $A_2$, $B_2$, $C_2$.  Then
we have a bijection between matrices $\gamma~\in~\Gamma$ and invariants
satisfying
$$(A_1,B_1,C_1) = (A_2,B_2,C_2) = 1,$$
$$A_1C_2 + B_1B_2 + C_1A_2 = 0,$$
$$A_1 \equiv A_2 \equiv B_1 \equiv B_2 \equiv 0 \, (3), C_1 \equiv C_2
\equiv 1\, (3).$$
We may now define the appropriate Kubota homomorphism
$\kappa: \Gamma \to \mu_3$.  Given $A_1$, $B_1$, $C_1$, $A_2$, $B_2$, $C_2$
as above, we factor (with $r_1~\equiv~r_2~\equiv 1~\,(3)$, $(C_1',C_2') = 1$)
$$B_1 = r_1B_1' \qquad B_2 = r_2B_2'$$
$$C_1 = r_1r_2C_1' \qquad C_2=r_1r_2C_2'.$$
The function
$$\kappa(g) =   \fracwithdelims(){B_1'}{C_1'}
                \fracwithdelims(){B_2'}{C_2'}
                {\fracwithdelims(){C_1'}{C_2'}}^{-1}
                \fracwithdelims(){A_1}{r_1}
                \fracwithdelims(){A_2}{r_2}$$
is independent of the choice of factorization and is the desired homomorphism
(cf.~\cite{BFH3, Theorem 2}).

Metaplectic forms on the cubic cover of $GL(3,\C)$ are just complex--valued
functions
on $GL(3,\C)/Z \cdot K$
which satisfy
$$f(\gamma \tau) = \kappa(\gamma)f(\tau) \text{ for all }
\gamma \in \Gamma, \tau \in GL(3,\C)/Z \cdot K \tag{**}$$
and which also satisfy a particular differential equation (not the
same equation as for forms on the double cover of $GL(2,\C)$, but
there is a global viewpoint which relates the two equations).  For the
rest of this paper, when we write ``satisfying condition (**),'' (or
condition (*)) we mean
satisfying condition (**) (or condition (*))
$and$ also satisfying that specified
differential equation.

\head
4. A certain metaplectic form
\endhead

In a fundamental paper, Bump and Hoffstein~\cite{BH} constructed an
Eisenstein series on the cubic cover of $GL(3,\C)$ which had the property
that its Whittaker--Fourier coefficients contained the Hecke L--series
of the cubic residue symbol.  As we noted in section 1, these L--series
are nearly the L--series of the elliptic curves we wish to study; indeed,
once
twisted by a $grossencharacter$, they become exactly the L--series we wish
to study.  Our goal in this section, then, is to review the theory
of metaplectic forms on $GL(3,\C)$ and the construction of
Bump and Hoffstein,
and then to modify this construction to obtain a form from which we
will be able to obtain information about the L--series of the curves
$E_D$.

We begin by setting notation, which will be constant throughout the rest
of the paper.
Let $Z$ denote the scalar, and $K$ the unitary matrices in $GL(3,\C)$,
as usual.  We now write $\Gamma^2$ for the copy of $SL(2,\O)$ embedded in
the upper left corner of $SL(3,\O)$, and $\Gamma^2_\infty$ for
the intersection
of $\Gamma^2$ with the group of unipotent upper
triangular matrices.  We continue to let
$\Gamma$ denote the discrete group defined in the previous
section, and write $\Gamma_P$ for those elements of $\Gamma$ which
have a bottom row of $(0~0~1)$ (this is the maximal parabolic subgroup
of $\Gamma$). If $\phi$ is a form
satisfying (**), then we have a Fourier expansion
$$\phi(\tau) = \sum_n \phi_{0,n}(\tau) +
\sum \Sb \gamma \in \Gamma^2_\infty \setminus \Gamma \\
         m,n; m \ne 0 \endSb
\phi_{m,n}(\gamma\tau).$$
The multiplicity-one theorem of Shalika asserts that
$$\phi_{m,n}\left(
\pmatrix y_1y_2 & & \\ & y_1 & \\ & & y_1 \endpmatrix
\right) = \dfrac{a_{m,n}}{\N(mn)}W(my_1,ny_2)$$
where $W$ is the appropriate ``Whittaker function'' (more on Whittaker
functions below).
What Bump and Hoffstein did was to construct a metaplectic
Eisenstein series $\phi$ which had the property that, for $mn \ne 0$,
$$a_{m,n} = N_{m,n} \times \sum \Sb a \in \O \\ a \equiv 1 \, (3) \endSb
\fracwithdelims(){m^2n}{a}_3
 \bold N (a)^{-s + \frac{1}{2}} \tag{***}$$
where $N_{m,n}$ is a finite sum, which increases in complexity with
the divisibility of $m$ and $n$ by prime powers.
(We will see below that since
$\phi$ is an Eisenstein series, it includes a (surpressed,
in this notation) parameter $s$ which is the parameter appearing in
the right hand side of (***).)

\demo{Remark}
In order to
obtain the particular L-series in (***) one must make a
different choice of cubic symbol than the one actually made in
\cite{BH}.  We will fix such a choice throughout.  In
particular, our symbol is conjugate to theirs.
\enddemo

In words, what Bump and Hoffstein did was construct the Eisenstein
series on $GL(3,\C)$ by inducing up from the cubic theta function on
$GL(2,\C)$.  (cf. the last section of this paper.  Their paper also contained
 several other fundamental results,
which we will not cover here.  The cubic theta function is
a function on the cubic cover of $GL(2,\C)$; see~\cite{Pa} for details).
More precisely, they studied the form
$$\phi(\tau) = \sum_{\gamma \in \Gamma_P \setminus \Gamma}
\kappa(\gamma) I_s(\gamma \tau)$$
where
$$I_s(\tau) = \left(y_1^2 y_2\right)^2 \theta_3 \left(\pmatrix
y_2 & x_2 \\ 0 & 1 \endpmatrix\right)$$
(Here, $y_i$ and $x_i$ are the usual Iwasawa
coordinates of $\tau$ mentioned in the last section; $\theta_3$ is
the $GL(2,\C)$ cubic theta function.  Note that $I_s(\tau)$ is invariant
under $\Gamma_P$ (since $(y_1^2y_2)$ is).)
We wish to modify the Bump-Hoffstein form in two ways.  First, we wish
to introduce the $grossencharacter$ into the Whittaker-Fourier coefficients
to construct a form which has the L--series of the desired elliptic curves
in its Whittaker-Fourier coefficients; second, we will (in the construction)
sum over a larger group, thus increasing the size of the discrete group
$\Gamma$  under which our form will be invariant.

\demo{Remark}
This last modification will make the machinery of our next section much
simpler than it would be otherwise.  In general, one would like to
make $\Gamma$ as large as possible, but there are in fact limitations on
the size of $\Gamma$; these are imposed by the requirement that a
Kubota homomorphism exist.
\enddemo

We begin by introducing the $grossencharacter$ (see~\cite{L} for
details of this process).  The first step
is to define a representation $\tp$ on the group
$$\left\{ \pmatrix A & B & \\ C & D & \\ & & \xi \endpmatrix \in
\K(3, \C) \quad \bigg| \quad
\xi = \frac{R}{|R|}, \text{ where } R \in \Cal
O, R \equiv 1, 2 \, (3) \right\} \tag{4.1}$$
such that
$$\tp\pmatrix A & B & \\ C & D & \\ & & \frac{R}{|R|}\endpmatrix
= \pm \frac{R}{|R|}$$
(the sign is positive precisely when $R \equiv 1 \, (3)$.)  Next, we
induce $\tp$ up to the full group $\K(3,\C)$ and let $\p$ be any
(fixed now, once and for all) finite dimensional subrepresentation of
this induced representation.

By abuse of notation, we will write $\p$ applied to a
scalar matrix $z$ times a
unitary matrix $k$.  In these instances, what we mean is to rewrite
the product $z \cdot k$ as the product of a positive real scalar
matrix times a unitary matrix (by multiplying $z$ by the
appropriate complex number of absolute value 1, and $k$ by the
conjugate of this number), and then to apply $\rho$ to $k$.

{}From Frobenius reciprocity, $\p$ restricted to
(4.1) contains at least one copy of the original
representation $\tilde{\p}$.  Further, since $\K(3,\C)$ is
compact, we can define a functional
$\T$  on $\p$ such that
$\T(\p(\tau)) = \tilde{\p}(\tau)$ if $\tau \in
(4.1)$, and $$\T(\p(\tau \cdot g)) = \T(\p(\tau) \p(g)) =
\tilde{\p}(\tau)\T(\p(g)) \text{ for } \tau \in (4.1).$$
For such
$\tau$ it is also true that $\T(\p(g \cdot \tau)) =
\tilde{\p}(\tau)\T(\p(g))$.

We may now finally define our ``I-function'' in terms of the
$\I_w'(\tau)$ introduced by Bump and Hoffstein.  For $g \in \GL(3,
\C)$, fix an Iwasawa  decomposition $g = \tau \cdot z \cdot
k$, where $z \in \Z(3, \C)$, $k \in \K(3, \C)$, and define
$$\I_w^\phi(g) = \I_w'(\tau) \T(\p(z \cdot k)). $$
That this is well-defined follows from properties of $\p$.
By abuse of notation, we will often
refer
to $\T(\p(\gamma))$ as
$\p(\gamma)$. We also note that $\I_w^\phi(\gamma g) = \I_w^\phi(g)$ for
$\gamma \in
\Gamma_P$.

We are now ready to define our main object of study in this
paper.  Let $\Gamma^{\Z}$ denote those elements of $\SL(3,
\Cal O_{\Q(\sqrt{-3})})$ which are congruent to an element
of $\SL(3, \Z)$ mod 3.  In addition, let $\Gamma_P^{Z}$
be those elements of $\Gamma^{\Z}$ which are of the
form $\pmatrix * & * & * \\ * & * & * \\ 0 & 0 & 1 \endpmatrix$.
We may now define the Eisenstein series

$$\phi(\tau,w) = \sum_{\gamma \in \Gamma_P^{\Z} \setminus
\Gamma^{\Z}} \kappa (\gamma) \I_w^\phi(\gamma \tau) $$
and the normalized Eisenstein series
$$\multline
\phi^*(\tau,w) =
\left(2\pi\right)^{1-3w} \\
\times \zeta\left(3w-\frac{4}{3}\right)
\zeta\left(3w-\frac{2}{3}\right) \Gamma \left(
\frac{3w}{2} - \frac{2}{3} \right) \Gamma \left( \frac{3w}{2} -
\frac{1}{3} \right) \sum_{\gamma \in \Gamma_P^{\Z} \setminus
\Gamma^{\Z}} \kappa(\gamma) \I_w^\phi(\gamma \tau)
\endmultline $$
where $\kappa(\gamma)$ the Kubota symbol defined and studied by
\cite{BFH} (recall that our cubic residue symbol is conjugate
to theirs).

For the rest of this
section, we fix the notation
$\phi^*(\tau,w)$ for this one particular Eisenstein series.

Recall the definition of the $m, n\nth$ Fourier
coefficient of $\phi^*(\tau,w)$,

$$\phi^*(\tau,w)_{m, n} = \int_{\C /3}\int_{\C /3}\int_{\C /3}
\phi_w^*\left(\pmatrix 1 &
\xi_2 & \xi_3 \\ & 1 & \xi_1 \\ & & 1 \endpmatrix \tau\right)
e(-m\xi_1 - n \xi_2) \, d\xi_1 \, d\xi_2 \, d\xi_3. \tag{4.2}$$
\rm (Note that the measures in (4.2) are complex Haar
measures).  This definition is actually quite natural, as these
Fourier coefficients will appear in the evaluation of certain
Rankin--Selberg convolutions.  The following theorem shows the relationship
of the Whittaker--Fourier coefficients of $\phi^*(\tau,w)_{m, n}$ to
to the L--series of the elliptic curves (0.1).

\proclaim{Theorem 3}
Suppose $mn \ne 0$.  The
$m,n\nth$ Fourier coefficient
$\phi^*(\tau,w)_{m, n}$
is equal to
$$ \T(w)_{m.n} \cdot \A(w)_{m,n} \cdot |mn|^{-2} \cdot
\W(w,my_1, ny_2) $$  where $\T(w)_{m,n}$ is the Dirichlet polynomial
$$\aligned
&\quad \sum_{\alpha \ge 0}\quad \sum_{\delta \le
\alpha}
3^{\left[\frac{\alpha-2\delta+1}{3}\right]}
\fracwithdelims(){3}{m'}^{\alpha-2\delta}
\left(3^\alpha\right)^{-1-3w} \\
 \cdot&\sum_{b_2 \text{ mod } 3^{\alpha-\delta}}
\fracwithdelims(){3^{\alpha-\delta}}{b_2}
e\left(n\frac{b_2}{3^{\alpha-\delta}}\right)
\sum_{C_2 \text{ mod } 3^\alpha}
\fracwithdelims(){3^\delta}{C_2}
e\left(m\frac{b_2^{-1}C_2}{3^\delta}\right) \\
+ \, 3& \sum_{\alpha \ge 0}
3^{\left[\frac{\alpha+1}{3}\right]}
\fracwithdelims(){3}{m'}^{\alpha}
\left(3^\alpha\right)^{-1-3w}
\sum_{b_2 \text{ mod } 3^{\alpha}}
\fracwithdelims(){3^{\alpha}}{b_2}
e\left(n\frac{b_2}{3^{\alpha}}\right) \\
+ \, 9&;
\endaligned $$
$\W(w,y_1,y_2)$ is the Whittaker
function
$$ \aligned &\pigam \\
 &\cdot \left(y_1^2y_2\right)^w y_2
\int_{\C}\int_{\C}
\rho \pmatrix -\xi_2Q^{-\frac{1}{2}} & y_2Q^{-\frac{1}{2}} & 0 \\
-y_2\xi_3P^{-\frac{1}{2}}Q^{-\frac{1}{2}} &
 -\xi_2\xi_3P^{-\frac{1}{2}}Q^{-\frac{1}{2}} &
y_1P^{-\frac{1}{2}}Q^{\frac{1}{2}} \\
y_1y_2P^{-\frac{1}{2}} & y_1\xi_2P^{-\frac{1}{2}}
  & \xi_3P^{-\frac{1}{2}}  \endpmatrix \\
& \qquad \qquad \qquad \times P^{\frac{1-3s}{2}}Q^{-1}K_{\frac{1}{3}}
\left(4 \pi y_2P^{\frac{1}{2}}Q^{-1}\right)
e\left(-\bar{\xi_2}\xi_3Q^{-1}\right)
e(- \xi_2) \, d\xi_2 \, d\xi_3; \endaligned \tag{4.3}$$
and $\A(w)_{m,n}$ is the product of
the L-series
$$\prod_{p; \, (p, \lambda^6 mn)=1, \, p \equiv 1 \, (3)}
\left(1-\fracwithdelims(){mn^2}{p}\frac{|p|}{p}\N
p^{1-\frac{3w}{2}}\right)^{-1}
\left(1-\left[\frac{|p|}{p}\right]^3\bold N
p^{3-\frac{9w}{2}}\right) $$
with a factor $$ \multline
\frac{1}{2}3^{\frac{5}{2}}  |mn^2|^{w-1} \\
\sum_{A' \equiv 1 \, (3); \, p|A'
\Rightarrow p | \lambda^6mn} \frac{|A'|}{A} |A'|^{1-3w}
\sum_{A'=a' \cdot d'; \, d' | \lambda^3ma'}
\tau\left(\frac{ma'}{d'}\right) \\
\times \sum_{b_1 \text{ mod } a'}
\fracwithdelims(){b_1}{a'}e\left(n\frac{b_1}{a'}\right)
\sum_{C_1 \text{ mod } d'}\fracwithdelims(){C_1}{d'}
e\left(m\frac{b_1^{-1}C_1}{d'}\right).
\endmultline $$

\endproclaim

\head
5. The Rankin--Selberg method
\endhead

A Rankin--Selberg convolution, in a general sense, is a machine for
constructing a Dirichlet series with coefficients which are built out of
the Whittaker--Fourier coefficients of automorphic forms, and for providing
analytic information about the behavior of that Dirichlet series.  This
information can then be translated, as we mentioned in section 1, into
information about the collective behavior of the coefficients of the
original forms.

The type of convolution we will discuss below is an example of a convolution
discovered independently by Asai~\cite{A} and Patterson~\cite{Pa}.  It has
one particularly nice feature: recall that the Whittaker--Fourier coefficients
of the form constructed in the previous section were indexed by integers
in a quadratic field.  The Asai--Patterson convolution produces a ``sieved''
Dirichlet series, made up of those coefficients with rational integral
indices.  For our application, this is precisely what we desire.

Recall that a $GL(3,\C)$ metaplectic form
automorphic with respect to $GL(3,\Cal O)$
has a Fourier expansion
$$ \phi(\tau) = \sum_{n \in \Cal O} \phi_{0,n}(\tau) +
\sum_{\gamma \in \Gamma^2_\infty(\C) \setminus \Gamma^2(\C)} \,
\sum \Sb m,n \in \Cal O \\ m \ne 0 \endSb
\phi_{m,n}(\gamma\tau). $$
and that we have
$$\phi_{m,n}(\tau) = \text{arithmetic part } \times
W_{m,n}(y_1,y_2) e(mx_1) e(nx_2)$$
where $e(\cdot)$ is an exponential function and $W_{m,n}(y_1,y_2)$
is Whittaker function.  If $mn~\ne~0$, then $W_{m,n}(y_1,y_2)$
is of rapid (at least exponential) decay in both $y_1$ and $y_2$ as
$y_1~\to~\infty$ and $y_2~\to~\infty$.  When $m$ or $n$ is zero,
then $W_{m,n}(y_1,y_2)$ is of polynomial growth/decay in $y_1$ or
$y_2$, respectively.  Note also that if $mn~\ne~0$, then
$W_{m,n}(y_1,y_2)=W(my_1,ny_2)$; that is, the non-degenerate
Whittaker functions are all the same.

If we define the functions $I_s(\tau) = \left
(y_1y_2\right)^s$ and
$$\zeta^*(s) =
\zeta(s)\pi^{-\frac{s}{2}}\Gamma\left(\frac{s}{2}\right),
$$then the maximal parabolic Eisenstein series
$$E(\tau,s) = \zeta^*(3s)
\sum_{\gamma \in SL(3,\Z)_P \setminus SL(3,\Z)}
I_s(\gamma \tau)$$ (again, $SL(3,\Z)_P$ denotes those elements
of $SL(3,\Z)$ with bottom row $(0\,\, 0\,\, 1)$) converges for
$re(s) > 1$, has a meromorphic continuation to the entire complex
plane with simple poles at $s=0$ and $s=1$, and has a functional
equation as $s \mapsto 1-s$.

When $\phi$ is a $GL(3,\C)$ cusp form ($\phi_{m,n}=0$ whenever
$mn=0$), then the convolution
$$R(s,\phi) = \int_{SL(3,\Z) \setminus (GL(3,\R)/Z \cdot K)} \phi(\tau)
E(\tau,s) d^h\tau$$
(where  $d^h\tau$ is the invariant Haar measure
$$d^h\tau = \frac{dx_1 dx_2 dx_3 dy_1 dy_2}{(y_1y_2)^3})$$
is convergent (since the non-degenerate
Whittaker functions are of rapid decay) and inherits the
analytic properties (in $s$) of $E(\tau,s)$.  On the other hand,
it is straightforward to verify that by unfolding the
Eisenstein series, and then unfolding the Fourier expansion of
$\phi$ (recall the inner sum over $GL(2,\C)$ matrices), we obtain
$$R(s,\phi) = \sum \Sb m,n \in \Z \\ mn\ne 0 \endSb
\frac{a_{m,n}}{|m^2n|^s} \int_0^\infty\int_0^\infty W(y_1, y_2)
y_1^{2s-3} y_2^{s-3} \, dy_1 dy_2.$$
We know that $R(s,\phi)$ has a pole at $s=1$, and so
if we show that the double Mellin transform of the Whittaker
function converges at $s=1$, then we have established that the
Dirichlet series
$$\sum \Sb m,n \in \Z \\ mn\ne 0 \endSb
\frac{a_{m,n}}{|m^2n|^s}$$ has a pole at $s=1$, and we may apply
the ``easy method'' of Section 2 to translate this information into
information about the coefficients $a_{m,n}$.  We may also use
the functional equation of $R(s,\phi)$ to apply Tauberian type
theorems to this Dirichlet series.

When $\phi$ is not a cusp form, the convolution above does not converge.
Nonetheless, it is possible to define a Rankin--Selberg convolution which
evaluates to the desired Dirichlet series, and to obtain (albeit with
much more work than in the cuspidal case) the functional equation and
meromorphic continuation of this convolution.  For $GL(2,\R)$, this
extension of the theory of Rankin--Selberg convolutions to non--cuspidal
forms was accomplished by Zagier~\cite{Z}.  The $GL(3,\C)$ convolution
was extended to non--cuspidal forms in~\cite{L2}, whose results we
summarize below.

Let $F_u$ denote the sum of the degenerate terms of the
Fourier expansion of $\phi$ (expanded with respect to the upper
parabolic subgroup), namely
$$F_u =\sum_{n \in \Cal O} \phi_{0,n}(\tau) +
\sum_{\gamma \in \Gamma^2_\infty(\C) \setminus \Gamma^2(\C)} \,
\sum_{m \in \Cal O, m \ne 0}
\phi_{m,0}(\gamma\tau). $$

We now may define the Rankin-Selberg convolution $R(\phi,w)$ by
$$R(\phi,s) =
\zeta^*(3s) \int_{SL(3,\R)_P \setminus (GL(3,\R)/Z \cdot K)} [\phi(\tau) -
F_u(\tau)] \left(y_1^2y_2\right)^s \, d^h\tau. $$
Note that this is analogous to unfolding the Eisenstein series in the
cuspidal convolution, and that it is still not difficult to evaluate this
convolution directly, and obtain the desired Dirichlet series.  The
meromorphic continuation and functional equation are given by the
following Theorem. (In the Theorem, we use the standard notation of
$\tilde{\phi}(\tau)~=~\phi(\iota(\tau))$.)

\proclaim{Theorem 4} \sl The Rankin-Selberg convolution converges
for $Re(w)$ sufficiently large, has a meromorphic continuation
to all $w$, and satisfies the functional equation
$$ R(\phi,w)=R(\tilde{\phi},1-w). $$
If the $0,0^{th}$ Fourier coefficient of $\phi$ is the polynomial
$$\phi_{0,0} = \sum_{i=1}^k y_1^{\alpha_i}y_2^{\beta_i}$$
then $R(\phi,w)$ may have poles only at the locations $0,
\frac{1}{3}, \frac{2}{3}, 1$, $\alpha_i - 1$, $2- \beta_i$, $1 -
\frac{\alpha_i}{2}$ and $\frac{\beta_i}{2}$.

Further, if the Fourier expansion of $\phi$ is given by (3.2)
where for $mn \ne 0$ we have
$$\phi_{m,n}(\tau) =
\frac{a_{m,n}}{|mn|^2} W(y_1, y_2)$$
then we may explicitly obtain
$$R(\phi, w) = \sum \Sb m,n \in \Z \\ mn \ne 0 \endSb
\frac{a_{m,n}}{|m^2n|^w} \int_0^\infty\int_0^\infty W(y_1, y_2)
y_1^{2w-3} y_2^{w-3} \, dy_1 dy_2.$$
\endproclaim

Applying this to our form of the last section, and checking
the decomposition (cf.~\cite{L}) of $R(\phi,s)$ into concrete integrals
for pole locations, we see that $R(\phi^*,s)$ appears to have a pole
of order $2$ at $s=1$.  The following Lemma verifies
that this pole comes from the
Dirichlet series part of the convolution and not from the double
Mellin transform of the Whittaker function.
\proclaim{Lemma}
Fix $w=1$.  Then the integral
$$\int_0^\infty \int_0^\infty \W (y_1, y_2) y_1^{2s-3} y_2^{s-3}
\, dy_1 dy_2$$
converges for $re(s) \ge 1$.
\endproclaim

We now wish to show that the pole comes from contributions from infinitely
many curves;  the problem is that
$$x^3 + y^3 = D \qquad  \text{and} \qquad x^3 + y^3 = m^3D$$
are the same curve.  We must show that as the sum over $L(E_{m^3D},1)$
is finite for each $D$.

\proclaim{Lemma}
Fix an cube--free integer $k$, and set $w=1$.
Assume
that $re(s) > \frac{42}{45}$.  Then
$$\sum \Sb m,n \in \Z^* \\ \text{the cube--free part of
$m^2n=k$} \endSb \frac{\T(w)_{m,n}\A(w)_{m,n}}{|m^2n|^s} \,\, < \,\, \infty$$
where ``the cube--free part of $m^2n=k$'' means that $m^2n$ is
equal to $k$ times a cube.
\endproclaim

We have thus established Theorem 1.  To obtain Theorem~2, we need a
convolution which produces the ``twisted'' Dirichlet series
$$\sum_{m,n}\frac{L(E_{m^2n},1)\chi(m^2n)}{|m^2n|^w}$$
where $\chi$ is an arbitrary multiplicative character.  For such convolutions,
the functional equation has not been established.  Nonetheless, the convergence
in a half--plane, and the pole locations have been established (cf.~\cite{L2}),
yielding Theorem~2.  This convolution is technically more difficult
because the once $\phi$ is twisted by a character, it is automorphic
under a smaller group.  The Eisenstein series in the convolution is thus
constructed from a smaller group (recall that we need to ``unfold'' this
Eisenstein series using the automorphy of the form in the convolution),
and has a more complicated functional equation.  This is precisely why
we summed over a larger group (in constructing $\phi^*$) than Bump and
Hoffstein did: so that we could use the maximal parabolic Eisenstein
series constructed from the full group $SL(3,\Z)$, and not from a
congruence subgroup.

The order of the pole of the Dirichlet series made up
from the $L(E_D,1)$ is quite important.  We have seen that it has at most a
double
pole.  The question of whether it actually has a double pole, or merely a
single
pole, depends on whether the double Mellin transform of the twisted Whittaker
function (4.3) vanishes.  If it does, the Dirichlet series has a single pole
(it is direct to verify that the single pole does indeed have a non-zero
residue).
If not, the Dirichlet series has a double pole (the Lemma above states that
this double transform converges, but says nothing about its actual value).
This problem is quite hard: indeed, one has almost no information about the
representation $\rho$ with which to work.  It might be possible to
answer this question directly (in the language above), but there is another
viewpoint which might provide a more insightful answer -- that of
representation theory, cf. the next Section.

\head
6. Future directions
\endhead

There are two main paths to generalize the work above.
One is to try to sharpen the analytic machinery in order to obtain
better average value results, or average value results for $L(E_D,1)$ where
$D$ ranges over smaller sets.  The other is to try to find Eisenstein
series which have other interesting L--series in their Whittaker--Fourier
coefficients; indeed, a nice long--term goal might be to try to find
the correspondence between L--series and the associated metaplectic
forms.

We start by considering this latter problem.
If we think of a metaplectic form as a form on $GL(r,\C)$ which satisfies
a certain transformation property (such as (*) or (**)), then we always
have a Fourier expansion (due to the
multiplicity--one formula of Shalika, which
asserts
that the space of Whittaker functions is one--dimensional for
$GL(r,\C)$).  It is not always possible, however, to compute explicitly
the arithmetic part of the Whittaker--Fourier coefficients, as it was in
section 4.  The reason for this is that if one defines metaplectic forms
directly as automorphic forms on the metaplectic group, in general there
is no uniqueness of Whittaker functions.  So to have any hope of
finding metaplectic forms with interesting L--series in their Fourier
coefficients, we would have have to prove the forms in question
had unique Whittaker functionals.  Conversely, if one could find forms
with unique Whittaker models, it would be interesting to compute their
Fourier coefficients.

Bump and Hoffstein~\cite{H} have conjectured that the Hecke L--series of
the $n\nth$ order residue symbol should appear in the Fourier coefficients
of the Eisenstein series obtained by inducing the theta function on the
$n-$cover of $GL(n-1,\C)$ up to the $n-$cover of
$GL(n,\C)$.  (The case $n=3$ is considered
in their paper~\cite{BH}.)  Bump and Lieman~\cite{BL} have shown that the
corresponding local representation (the representation
obtained by inducing the exceptional
representation on the $n-$cover of $GL(n-1,\F)$ up to the $n-$cover of
$GL(n,\F)$, where $\F$ is a local field) does
indeed have unique Whittaker models.  This result
provides evidence that the conjecture of Bump and Hoffstein is provable
(indeed, this conjecture is the target of an ongoing investigation of Bump
and Lieman).  (Further, the viewpoint of this approach (representation theory)
is the correct one to determine whether the twisted Whittaker function,
discussed at the end of Section 5, does indeed vanish.)
Once this conjecture is proved, one could obtain the L--series
of biquadratic twists of elliptic curves as Whittaker--Fourier coefficients
of a metaplectic form on the $4-$cover of $GL(4,\C)$ by introducing a
$grossencharacter$ as above.

The other main avenue for future research is
to develop improved analytic machinery.
Farmer, Kumanduri and Lieman are working to obtain the functional equation
for the twisted $GL(3,\C)$ Asai--Patterson convolution discussed at the end
of Section 5; among other applications, this would allow one
to obtain average values of of $L(E_D,1)$ where $D$
ranges over rational integers in
a fixed congugacy class modulo any prime $p \ne 3$.  Farmer,
Hoffstein and Lieman are also working on obtaining the
functional equation for the L--series of $GL(3,\C)$ forms twisted
by a character.  This will allow the computation of the
average values of the L--series of $E_D$ where $D$ ranges not
over rational integers, but over all cube--free
integers in $\Cal O$.


\newpage
\Refs
\widestnumber\key{GHP-S}

\ref\key A
\by T. Asai
\paper On certain Dirichlet series associated with Hilbert modular forms
and Rankin's method
\jour Math. Ann.
\vol 226
\year 1977
\pages 81--94
\endref

\ref\key BFH1
\by D. Bump, S. Friedberg,and J. Hoffstein
\paper Eisenstein series on the metaplectic group and nonvanishing
theorems for automorphic L--functions and their derivatives
\jour Annals of Math
\vol 131
\pages 53--127
\year 1990
\endref

\ref \key BFH2
\by D. Bump, S. Friedberg, and J. Hoffstein
\paper Nonvanishing theorems for L--functions of modular forms and their
derivatives
\jour Invent. math.
\vol 102
\pages 543--618
\yr 1990
\endref

\ref\key BFH3
\by D. Bump, S. Friedberg, and J. Hoffstein
\paper Some cubic exponential sums
\jour preprint
\endref

\ref \key BH
\by D. Bump and J. Hoffstein
\paper Cubic metaplectic forms on GL(3)
\jour Invent. math
\vol 84
\yr 1986
\pages 481--505
\endref

\ref \key BL
\by D. Bump and D. Lieman
\paper Uniqueness of Whittaker functionals on the metaplectic group
\jour Duke Math. J., to appear
\endref

\ref\key BZ
\paper Representations of the group GL(n,F) where F is a non-archimedean
local field
\by I. N. Bernshtein and A. V. Zelevinskii
\jour Russian Math. Surveys
\vol 31:3
\yr 1976
\pages 1--68
\endref

\ref\key BZ2
\paper Induced representations of reductive p-adic groups
\by I. N. Bernshtein and A. V. Zelevinskii
\jour Ann. Scient. Ec. Norm. Sup.
\vol 10
\yr 1977
\pages 441--472
\endref

\ref\key D
\by L. Dickson
\book History of the Theory of Numbers II
\yr 1934
\publ G.E. Stechert \& Co.
\endref

\ref\key FL
\by Y. Flicker
\paper Twisted tensors
and Euler products
\jour Bull. Soc.
math. France
\vol 116
\pages 395--313
\yr 1988
\endref

\ref\key FR
\by S. Friedberg
\paper A global approach to the Rankin-Selberg
convolution for $GL(3, \Z)$
\jour Trans. AMS
\vol 300
\pages 159--178
\yr 1987
\endref

\ref \key GHP-S
\by S. Gelbart, R. Howe, and I. Piatetski-Shapiro
\paper Uniqueness and existence of Whittaker models for the metaplectic group
\jour Israel Journal of Mathematics
\vol 34
\yr 1979
\pages 21--37
\endref

\ref\key GHP
\by D. Goldfeld, J. Hoffstein, and S. Patterson
\paper On automorphic functions of half-integral weight with
applications to elliptic curves
\jour in ``Number theory related to Fermat's Last Theorem,'' N.
Koblitz, editor
\pages 153--194
\year 1982
\endref

\ref\key GV
\by D. Goldfeld and C. Viola
\paper Some conjectures on elliptic curves over cyclotomic fields
\jour Trans. AMS
\vol 276
\year 1983
\pages 511--515
\endref

\ref \key H
\by J. Hoffstein
\paper Eisenstein series and theta functions on the metaplectic group
\jour preprint
\endref

\ref\key IR
\by K. Ireland and M. Rosen
\book A Classical Introduction to Modern Number Theory
\publ Springer-Verlag, Graduate Texts in Math.
\vol 84
\yr 1982
\endref

\ref \key KP
\by D. Kazhdan and S. Patterson
\paper Metaplectic Forms
\jour Publ. Math. IHES
\vol 59
\yr 1984
\pages 35--142
\endref

\ref \key L
\by D. Lieman
\paper
Nonvanishing of L--series associated to cubic twists of elliptic curves
\jour Annals of Mathematics (to appear)
\endref

\ref \key L2
\by D. Lieman
\paper
The GL(3) Rankin--Selberg convolution for functions
not of rapid decay
\jour
Duke Mathematical Journal
\vol 69
\yr 1993
\pages 219--242
\endref

\ref\key N
\by J. Nekovar
\paper
Class numbers of quadratic fields and Shimura's correspondence
\jour Math Ann.
\vol 287
\yr 1990
\pages 577--594
\endref

\ref\key Pa
\by S. Patterson
\paper On Dirichlet series associated with cubic Gauss sums
\jour J. Reine Angew. Math.
\vol 303
\year 1978
\pages 102--125
\endref

\ref\key Pr
\by I. Proskurin
\paper
Automorphic functions and the homomorphism of Bass-Milnor-Serre~I, II (in
Russian)
\jour
Zapiski Naucnik Seminarov LOMI
\vol 129
\yr 1983
\pages 85--126 and 127--163
\endref

\ref\key Z
\by D. Zagier
\paper The Rankin-Selberg method for automorphic functions
which are not of rapid decay
\jour J. Fac. Sci. Tokyo Univ
\vol 28
\pages 416--437
\year 1981
\endref

\ref\key ZK
\by D. Zagier and S. Kramarz
\paper Numerical
Investigations related to the L-series of certain elliptic
curves
\jour Journal of the Indian Math Soc.
\vol 52
\pages 51--69
\year 1987
\endref
\endRefs
\enddocument
\end